\documentclass{amsart}

\newtheorem{theorem}{Theorem}
\newcommand{\bt}{\begin{theorem}}
\newcommand{\et}{\end{theorem}}
\newtheorem{lemma}{Lemma}
\newcommand{\bl}{\begin{lemma}}
\newcommand{\el}{\end{lemma}}
\newtheorem{corollary}{Corollary}
\newcommand{\bc}{\begin{corollary}}
\newcommand{\ec}{\end{corollary}}
\newtheorem{problem}{Problem}
\newcommand{\bprob}{\begin{problem}}
\newcommand{\eprob}{\end{problem}}
\newcommand{\beq}{\begin{equation}}
\newcommand{\eeq}{\end{equation}}
\newcommand{\benum}{\begin{enumerate}}
\newcommand{\eenum}{\end{enumerate}}

\newcommand{\Z}{\ensuremath{\mathbf Z}}

\newcommand{\R}{\ensuremath{\mathbf R}}

\newcommand{\mcc}{\ensuremath{ \mathcal C}}
\newcommand{\mcd}{\ensuremath{ \mathcal D}}

\newcommand{\mcg}{\ensuremath{ \mathcal G}}
\newcommand{\mch}{\ensuremath{ \mathcal H}}

\DeclareMathOperator{\id}{id}

\newcommand{\bmat}{\left(\begin{matrix}}
\newcommand{\emat}{\end{matrix}\right)}

\DeclareMathOperator{\qqand}{\qquad\text{and}\qquad}

\title[MSTD sets]{Problems in additive number theory, V: \\
Affinely inequivalent MSTD sets}
\author{Melvyn B. Nathanson}
\address{Department of Mathematics\\
Lehman College (CUNY)\\
Bronx, NY 10468} 
\email{melvyn.nathanson@lehman.cuny.edu}


\date{\today}

\begin{document}
\maketitle

\begin{abstract}
An MSTD set is a finite set of integers with more sums than differences.  
It is proved that, for infinitely many positive integers $k$, there are infinitely  
many affinely inequivalent MSTD sets of cardinality $k$. 
There are several related open problems.   
\end{abstract}

\section{Sums and differences}
In mathematics, simple calculations often suggest hard problems.  
This is certainly true in number theory.  
Here is an example:
\[
3+2 = 2 + 3 
\]
but
\[
3-2 \neq 2-3.
\]
This leads to the following question.  Let $A$ be a set of integers,  a set of real numbers, 
or, more generally, a subset of an additive abelian group \mcg. 
We denote the cardinality of the set $A$ by $|A|$.  
Define the \emph{sumset} 
\[
A+A = \{a+a': a,a' \in A\}
\]
and  the \emph{difference set} 
\[
A - A = \{a - a': a,a' \in A\}.
\]
For all $a, a' \in \mcg$ with $a \neq a'$, we have  $a+a' = a'+a$ because \mcg\ is abelian.  
However, $a-a' \neq a'-a$ if \mcg\ is a group, such as \R\ or \Z, 
with the property that $2x = 0$ if and only if $x=0$.  
It is reasonable to ask:  In such groups, does every finite set  
have the property that the number of sums does not exceed the number of differences?   
Equivalently, is $|A+A| \leq |A-A|$ for every finite subset $A$ of \mcg? 

The answer is ``no.''  A set with more sums than differences is called an \emph{MSTD set}.  

As expected, most finite sets $A$ of integers do satisfy $|A+A| < |A-A|$ 
(cf. Hegarty and Miller~\cite{hega-mill09} and Martin and O'Bryant~\cite{mart-obry07a}).  
For example, if 
\[
A = \{0,2,3\}
\]
then 
\[
A+A = \{0,2,3,4,5,6\}
\]
\[
A-A = \{-3,-2,-1,0,1,2,3\}
\]
and 
\[
|A+A| = 6 < 7 =  |A-A|.
\]
It is also easy to construct finite sets $A$ for which the number of sums 
equals the number of differences.
For example, if $A$ is an arithmetic progression of length $k$  in a torsion-free abelian group, 
that is, a set of the form 
\beq  \label{PANT-V:AP}
A = \{a_0 + id: i = 0,1,2,\ldots, k-1\}
\eeq
for some $d \neq 0$, then the number of sums equals the number of differences: 
\[
A+A = \{a_0 + id: i = 0,1,2,\ldots, 2k-2\}
\]
\[
A-A = \{a_0 + id: i = -(k-1), -(k-2),\ldots, -1, 0,1,\ldots, k-2, k-1\}
\]
and 
\[
 |A+A| = |A-A| = 2k-1.
\]

In an abelian group \mcg, the set $A$ is \emph{symmetric} if there exists an element $w \in \mcg$ 
such that $a \in A$ if and only if $w-a \in A$.  
For example, the arithmetic progression~\eqref{PANT-V:AP} 
is symmetric with respect to $w = 2a_0 + (k-1)d$.  
We can prove that every finite symmetric set has the same number of sums and differences.  
More generally, for $0 \leq j \leq h$, consider the \emph{sum-difference set}
\[
(h-j)A - jA = 
\left\{ \sum_{i=1}^{h-j} a_i -  \sum_{i=h-j+1}^h  a_i : a_i \in A \text{ for } i=1,\ldots, h \right\}.  
\]
For $h=2$ and $j=0$, this is the sumset $A+A$.  
For $h=2$ and $j=1$, this is the difference set $A-A$.

\bl           \label{PANT-V:lemma:shortSumDifference}
Let $A$ be a nonempty finite set of real numbers with $|A| = k$.  
For $j \in \{ 0,1,2,\ldots, h\}$, there is the sum-difference inequality 
\[
 |(h-j)A - jA| \geq h(k-1)+1.
\]
Moreover, 
\[
 |(h-j)A - jA| = h(k-1)+1
\]
if and only if $A$ is an arithmetic progression.
\el

\begin{proof}
If $A$ is a set of $k$  real numbers, then $|hA| \geq h(k-1)+1$. 
Moreover, $|hA| = h(k-1)+1$ if and only if $A$ is an arithmetic progression 
(Nathanson~\cite[Theorem 1.6]{nath96bb}).  

For every number $t$, the translated set  $A' = A - t$ satisfies 
\[
(h-j)A' - jA' = (h-j)A - jA - (h-2j)t
\]
and so 
\[
|(h-j)A' - jA'|  = |(h-j)A - jA|.
\]
Thus, after translating by $t = \min(A)$, we can assume that $0 = \min(A)$.  
In this case, we have  
\[
(h-j)A \cup (-jA) \subseteq (h-j)A-jA.
\]
Because $(h-j)A$ is a set of nonnegative numbers and $-jA$ is a set of nonpositive 
numbers, we have 
\[
(h-j)A \cap (-jA) = \{0\}
\]
and so
\begin{align*}
| (h-j)A-jA | & \geq |(h-j)A| + |-jA| - 1 \\ 
& \geq ( (h-j)(k-1)+1 ) + (j(k-1)+1) - 1 \\ 
& = h(k-1)+1.
\end{align*}
Moreover, $| (h-j)A-jA | = h(k-1)+1$ if and only if 
both $|(h-j)A | = (h-j)(k-1)+1$ and $|-jA |= j(k-1)+1$,   
or, equivalently,  if and only if $A$ is an arithmetic progression.  
This completes the proof.  
\end{proof}

\bt            \label{PANT-V:theorem:symmetricSD}
Let $A$ be a nonempty finite subset of an abelian group \mcg.  
If $A$ is symmetric, then 
\beq                                         \label{PANT-V:symmetricSD}
 |(h-j)A - jA| = |hA| 
\eeq
for all integers $j \in \{ 0,1,2,\ldots, h\}$.  
In particular, for $h = 2$ and $j=1$,
\[
|A-A| = |A+A|.
\]
Thus, symmetric sets have equal numbers of sums and differences.   
\et

Note that the nonsymmetric set 
\[
A = \{0,1,3,4,5,8\}
\]
satisfies 
\[
A+A = [0,16] \setminus \{14,15\} 
\qqand
A-A = [-8,8] \setminus \{\pm 6\} 
\]
and so
\[
|A+A | = |A-A | = 15.
\]
This example, due to Marica~\cite{mari69}, shows that there also exist non-symmetric 
sets of integers with equal numbers of sums and differences.

\begin{proof}
If $j=0$, then $(h-j)A - jA = hA$.  If $j = h$, then $(h-j)A - jA = -hA$.  
Equation~\eqref{PANT-V:symmetricSD} holds in both cases. 
Thus, we can assume that $1 \leq j \leq h-1$.  

Let $A$ be a symmetric subset with respect to $w \in \mcg$. 
Thus, $a \in A$ if and only if $w-a \in A$.  
For every integer $j$, define the function $f_j:\mcg \rightarrow \mcg$ by $f_j(x) = x + jw$.   
For all $j,{\ell} \in \Z$ we have  $f_j f_{\ell} = f_{j+{\ell}}$.  
In particular, $f_j f_{-j} = f_0 = \id$ and $f_j$ is a bijection.

Let $x = \sum_{i=1}^h a_i  \in hA$, and  let  $a'_i = w-a_i \in A$ for $i=1,\ldots, h$.    
If $1 \leq i \leq j \leq h$, then  
\begin{align*}
f_{-j}(x) 
& = \left( \sum_{i=1}^h a_i   \right) - jw \\ 
& =  \sum_{i=1}^{h-j}  a_i  - \sum_{i=h-j+1}^h( w - a_i ) \\
& =    \sum_{i=1}^{h-j}   a_i  - \sum_{i=h-j+1}^h a'_i \\ 
& \in (h-j)A - jA
\end{align*}
and so
\[
|hA| \leq |  (h-j)A - jA|.
\]

Let $y =  \sum_{i=1}^{h-j}  a_i  -\sum_{i=h-j+1}^h  a_i \in  (h-j)A - jA$.
For $h-j+1 \leq i  \leq h$, let  $a'_i = w-a_i \in A$.  
Then 
\begin{align*}
f_j(y) 
& = \left(  \sum_{i=1}^{h-j}  a_i  -\sum_{i=h-j+1}^h  a_i \right) + jw \\ 
& = \sum_{i=1}^{h-j}  a_i  + \sum_{i=h-j+1}^h  (w - a_i) \\ 
& = \sum_{i=1}^{h-j}  a_i  + \sum_{i=h-j+1}^h a'_i \\ 
& \in hA
\end{align*}
and so
\[
|(h-j)A - jA| \leq |hA|.
\]
Therefore, $|(h-j)A - jA| = |hA|$ and the proof is complete.     
\end{proof}

Let $A$ be a nonempty set of integers.  
We denote by $\gcd(A)$ the greatest common divisor of the integers in $A$.   
For real numbers $u$ and $v$, we define the \emph{interval of integers} 
$[u,v] = \{ n \in \Z: u \leq n \leq v \}$.  If $u_1,v_1,u_2,v_2$ are integers, 
then $[u_1, v_1] + [u_2, v_2] = [u_1+u_2, v_1+v_2]$.  

\bt
Let $A$ be a finite set of  nonnegative integers with $|A| \geq 2$ such that 
$0 \in A$ and $\gcd(A) = 1$.  Let $a^* = \max(A)$.  
There exist  integers $h_1$, $C$, and $D$ 
and sets of integers $\mcc^* \subseteq [0,C+D-1]$ and $\mcd^* \subseteq [0,C+D-1]$ 
such that, if $h \geq 2h_1$, then the sum-difference set 
has the structure 
\[
ja^* + (h-j)A - jA   = \mcc^* \cup [C+D,  ha^* - (C+D)] \cup (ha^* - \mcd^* )
\]
for all  integers $j$ in the interval $[h_1, h-h_1]$.
Moreover, 
\[
|(h-j)A - jA |= |(h-j')A - j' A|
\]
for all integers $j, j' \in [h_1, h-h_1]$.
\et

\begin{proof}
Because $A \subseteq [0, a^*]$, we have 
$hA  \subseteq [0, ha^*]$ for all nonnegative integers $h$.  
By a fundamental theorem of additive number theory
(Nathanson~\cite{nath72e,nath96bb}), 
there exists a positive integer $h_0 = h_0(A)$ 
and there exist nonnegative integers $C$ and $D$ 
and sets of integers $\mcc \subseteq [0,C-2]$ and $\mcd \subseteq [0,D-2]$ 
such that, for all  $h \geq h_0$, the sumset $hA$ has the rigid structure 
\beq             \label{PANT-V:Fundamental}
hA = \mcc \cup [C, ha^* - D] \cup \left( ha^* - \mcd \right).  
\eeq
Let 
\beq                    \label{PANT-V:h1}
h_1 = h_1(A) = \max\left( h_0, \frac{2C+D}{a^*},  \frac{C+2D}{a^*}  \right).
\eeq
Let $h \geq 2h_1$.  If $j \in [h_1, h-h_1]$, then 
\[
j \geq h_1 \qqand h-j \geq h_1.  
\] 
Let $r = h-j$.  Applying the structure theorem~\eqref{PANT-V:Fundamental}, 
we obtain the  sumsets  
\[
rA = \mcc \cup [C, ra^* - D] \cup \left( ra^* - \mcd \right)  
\]
and 
\[
jA = \mcc \cup [C, ja^* - D] \cup \left( ja^* - \mcd \right).   
\]
Rearranging the identity for $jA$ gives 
\[
ja^* - jA = \mcd \cup [D, ja^* - C] \cup  \left( ja^* - \mcc \right). 
\]
We have 
\begin{align*}
[C+D,  ha^* - (C+D)] & =  [C, ra^* - D]  + [D, ja^* - C] \\
& \subseteq  rA + (ja^* - jA). 
\end{align*}

It follows from~\eqref{PANT-V:h1} that 
\begin{align*}
\min(ja^* - \mcc) & \geq ja^* - (C-2) \\
& > ja^* - C \\ 
& \geq   h_1a^* - C \\
& \geq (2C+D) - C = C+D.  
\end{align*}
Similarly, 
\[
\min(ra^* - \mcd) > ra^* - D \geq  C+D.
\]
These lower bounds imply that for 
\[
n \in [0,C+D - 1]  \qqand j \in [h_1, h-h_1]
\]
we have $n \in rA + (ja^* - jA) $ if and only if 
\[
n \in (\mcc  + \mcd) \cup (\mcc + [D, ja^* - C] ) \cup  (\mcd +  [C, ra^* - D] )
\]
 if and only if 
\[
n \in (\mcc  + \mcd) \cup (\mcc + [D, C+D] ) \cup  (\mcd +  [C, C+D] ).
\]
Therefore, 
\begin{align*}
\mcc^* & = [0,C+D-1] \cap \left( (\mcc  + \mcd) \cup (\mcc + [D, C+D] ) \cup  (\mcd +  [C, C+D] ) \right) \\
& =  [0, C+D-1]   \cap  (rA + (ja^* - jA)) 
\end{align*}
for all $j \in [h_1, h-h_1]$.    
Similarly, there exists a set $\mcd^* \subseteq [0, C+D-1]$ such that 
\[
ha^* - \mcd^*  =  [ha^* - (C+D) + 1), ha^*]   \cap  (rA + (ja^* - jA))
\]
for all $j \in [h_1, h-h_1]$.  
Therefore,
\begin{align*}
ja^* + (h-j)A - jA  & = (rA + (ja^* - jA)) \\ 
& = \mcc^* \cup [C+D,  ha^* - (C+D)] \cup (ha^* - \mcd^* )
\end{align*}
for all $j \in [h_1, h-h_1]$.  
This completes the proof.  
\end{proof}

\bprob
Let $A$ be a set of $k$ integers.  
For $j = 0,1,\ldots, h$, let 
\[
f_{A,h}(j) =  |(h-j)A - jA|.
\]
Is 
\[
\max(f_{A,h}(j):j=0,1,\ldots, h) = f_{A,h}\left(  \left[ \frac{h}{2} \right]\right) ?
\]
Is the function $f_{A,h}(j)$ unimodal?
\eprob

Although the conjecture that a finite set of integers 
has no more sums than differences is reasonable, 
the conjecture is false.  Here are three counterexamples.  The set 
\[
A = \{0, 2, 3, 4, 7, 11, 12, 14\}
\]
with $|A|  = 8$ and with sumset  
\[
A+A = [0,28] \setminus \{1,20,27\}
\]
and difference set 
\[
A - A = [-14,14] \setminus \{  6,-6, 13, -13 \}
\]
satisfies 
\[
|A+A| = 26 > 25 = |A -A|.
\]
Note that $A =  \{0, 2, 3,  7, 11, 12, 14\} \cup \{4\}$, where the set $\{0, 2, 3,  7, 11, 12, 14\}$ 
is symmetric.  This observation is exploited in Nathanson~\cite{nath07yzx}.  

The set 
\[
B = \{ 0, 1, 2, 4, 7, 8, 12, 14, 15 \}
\]
with $|B| = 9$ and with sumset  
\[
B+B = [0,30] \setminus \{ 25 \}
\]
and difference set 
\[
B - B = [-15,15] \setminus \{  9, -9 \}
\]
satisfies 
\[
|B+B| = 30 > 29 = |B -B|.
\]

The set 
\[
C = \{ 0, 1, 2, 4, 5, 9, 12, 13, 14, 16, 17, 21, 24, 25, 26, 28, 29 \}
\]
with $|C| = 17$ and with sumset  
\[
C+C = [0,58] 
\]
and difference set 
\[
C - C = [-29, 29] \setminus \{ \pm 6, \pm 18  \}
\]
satisfies 
\[
|C+C| = 59 > 55 = |C -C|.
\]
Set $B$ appears in Marica~\cite{mari69} and set $C$ 
in Freiman and Pigaev~\cite{frei-piga73}.

An \emph{MSTD set} in an abelian group \mcg\ is a finite set that  
has more sums than differences.  
MSTD sets of integers have been extensively investigated in recent years, 
but they are still mysterious and many open problems remain.  MSTD sets of real numbers and MSTD sets in arbitrary abelian groups have also been studied.  
In this paper we consider only MSTD sets contained in the additive groups \Z\ 
and  \R.  
There are constructions of various  infinite families of MSTD sets of integers 
(e.g. Hegarty~\cite{hega07a},    Miller, Orosz, and Scheinerman~\cite{mill-oros-sche10}, and  Nathanson~\cite{nath07xyz}), 
but there is no complete classification.

\bprob
A fundamental problem is to classify the possible structures of MSTD sets of integers 
and of real numbers.   
\eprob

Let \mcg\  denote \R\ or \Z.
For all $\lambda, \mu \in \mcg$ with $\lambda \neq 0$,  we define the \emph{affine map} 
$f:\mcg \rightarrow \mcg$ by 
\[
f(x) = \lambda x + \mu.
\]
An affine map is one-to-one.   
Subsets $A$ and $B$ of \mcg\ are \emph{affinely equivalent} if there exists an 
affine map $f:A \rightarrow B$ or $f:B \rightarrow A$ that is a bijection.

Let $k \geq 2$ and let $A = \{a_0,a_1,\ldots, a_{k-1} \}$ be a set of integers 
such that 
\[
a_0 < a_1 < \cdots < a_{k-1}.  
\]
Let 
\[
d = \gcd(\{a_i-a_0:i=1,\ldots, k-1\})
\]
and
\[
a'_i = \frac{a_i-a_0}{d}
\]
for $i=0,1,\ldots, k-1$.  
Let $A' =\{a'_0,a'_1,\ldots, a'_{k-1} \}$.  We have  
\[
0 = a'_0 < a'_1 < \cdots < a'_{k-1}.  
\]
Note that 
\[
\min(A') = 0 \qqand \gcd(A') = 1.
\]
We call $A'$ the \emph{normal form} of $A$.

Consider the affine map  $f(x) = dx+a_0$.  We have 
\[
A = \{da'_i + a_0: i=0,1,\ldots, k-1\} = \{f(a'_i): i=0,1,\ldots, k-1\} = f(A')
\]
and so $f:A' \rightarrow A$ is a bijection 
and the sets $A$ and $A'$ are affinely equivalent.  

A property of a set is an \emph{affine invariant} if, for all affinely equivalent sets 
$A$ and $B$, the set $A$ has the property if and only if 
the set $B$ has the property.  

The property of being an MSTD set is an affine invariant.  
Let $f$ be an affine map on \mcg.  
For all $a_{i_1}, a_{i_2 }, a_{i_3 }, a_{i_4 }  \in \mcg$,  the following statements are equivalent:
\[
a_{i_1}  -a_{i_2 }  = a_{i_3 } - a_{i_4 }  
\]
\[
a_{i_1} + a_{i_4 }   =  a_{i_2 }  + a_{i_3 } 
\]
\[
f(a_{i_1}) + f(a_{i_4}) = f(a_{i_2}) + f(a_{i_3}) 
\]
\[
f(a_{i_1}) - f(a_{i_2}) = f(a_{i_3}) - f(a_{i_4}).
\]
This implies that if $A$ is an MSTD set, 
then $B = f(A)$ is an MSTD set for every affine map $f$.
Thus, to classify MSTD sets of real numbers or of integers, it suffices to classify them 
up to affine maps.  

In the group of integers, 
Hegarty~\cite{hega07a}  proved that that there exists no MSTD set 
of cardinality less than 8, 
and that every MSTD set of cardinality 8 is affinely equivalent to the set 
$\{ 0,2,3,4,7,11,12,14 \}$. 

Let $\mch(k,n)$ denote the number of affinely inequivalent  
MSTD sets of integers of cardinality $k$ contained 
in the interval $[0,n]$.
Thus, Hegarty proved that $\mch(k,n) = 0$ for $k \leq 7$ and all 
positive integers $n$, that $\mch(8,n) = 0$ for $n \leq 13$, 
and that  $\mch(8,n) = 1$ for $n \geq 14$.

\bprob
Why does there exist no MSTD set of integers of size 7?  
\eprob

\bprob
Let $k \geq 9$.  Compute $\mch(k,n)$.  Describe the asymptotic growth of $\mch(k,n)$ 
as $n \rightarrow \infty$.
\eprob

\bprob
For fixed $n$, describe the behavior of $\mch(k,n)$ as a function of $k$.  
For example, is $\mch(k,n)$ a unimodal function of $k$?  
Note that $\mch(k,n) = 0$ for $k > n$.  
\eprob

For fixed $k$, the function $\mch(k,n)$ is a monotonically increasing function of $n$.
Denoting by $\mch(k)$  the number of affinely inequivalent  
MSTD sets of cardinality $k$, we have 
\[
\mch(k) = \lim_{n \rightarrow \infty} \mch(k,n).
\]
Thus, $\mch(k)  = \infty$ if there exist infinitely many affinely inequivalent 
MSTD sets of integers of cardinality $k$.

For every finite set $A$ of integers, define
\[
\Delta(A) = |A - A| - |A+A|.
\]
The set $A$ is an MSTD set if and only if $\Delta(A) < 0$.

\bl              \label{PANT-V:lemma:Delta}
Let $A = \{a_0,a_1,\ldots, a_{k-1} \}$ be a set of $k$ integers with 
\[
0 = a_0 < a_1 < \cdots < a_{k-1}. 
\]
If $a_k$ is an integer  such that 
\[
 2a_{k-1}  <  a_k 
\]
and if
\[
A' = A \cup \{ a_k\} 
\]
then 
\[
\Delta(A')  - \Delta(A) = k -1.
\]
\el

\begin{proof}
We have 
\[
A'+A' = (A+A) \cup \{a_k+a_i:i=0,1,\ldots, k\}.
\]
Because 
\[
\max(A+A) = 2a_{k-1} < a_k < a_k+a_1 < \cdots < a_k + a_{k-1}  <  2a_k
\]
we have 
\[
|A'+A'| = |A+A| + k+1.
\]

Similarly, 
\[
A'-A' = (A-A) \cup \{ \pm(a_k - a_i) :i=0,1,\ldots, k-1 \}.
\]
Because 
\[
\max(A-A) = a_{k-1} < a_k - a_{k-1} < a_k - a_{k-2} < \cdots < a_k - a_1 < a_k
\]
and
\[
\min(A-A) = -a_{k-1} > -a_k + a_{k-1} >  \cdots >   - a_k + a_1 > - a_k
\]
we have 
\[
|A'-A'| = |A-A| + 2k.
\]
Therefore,
\begin{align*}
 \Delta(A') & = |A'-A'| - |A'+A'| \\
 & =  (|A-A| + 2k) - ( |A+A| + k+1) \\
&  =  \Delta(A) + k-1.
\end{align*}
This completes the proof. 
\end{proof}

\bl              \label{PANT-V:lemma:Delta-k}
Let $B$ be an MSTD set  of  integers with 
\[
|B+B| \geq |B - B| + |B|.  
\]
There exist infinitely many affinely inequaivalent MSTD sets of integers 
of cardinality $|B|+1$, that is, 
$\mch( |B|+1) = \infty$.
\el

\begin{proof}
Let $|B| = {\ell}$.  
Translating the set $B$ by $\min(B)$, we can assume that $0 = \min(B)$.  
Let $b_{{\ell}-1}= \max(B)$.  
The inequality 
\[
|B+B| \geq |B-B| + |B| 
\]
is equivalent to  
\[
\Delta(B) \leq - {\ell}.  
\]
For every integer $b_{\ell} > 2 b_{{\ell}-1}$ and $B' = B \cup \{ b_{\ell} \}$, 
Lemma~\ref{PANT-V:lemma:Delta} implies that  
\[
\Delta(B') = \Delta(B) + {\ell} - 1 \leq -1
\]
and so 
\[
|B' - B'| < |B' + B'|.
\]
Therefore, $B'$ is an MSTD set of integers of cardinality ${\ell}+1$.  
If $b'_{\ell} > b_{\ell} > 2b_{{\ell}-1}$, then the sets $B \cup \{ b_{\ell} \}$ 
and $B \cup \{ b'_{\ell} \}$ are affinely inequivalent, 
and so $\mch({\ell}+1) = \infty$.
\end{proof}

\bl              \label{PANT-V:lemma:Delta-m}
Let $A$ be a nonempty finite set of nonnegative integers with $a^* = \max(A)$.
Let $m$ be a positive integer with 
\[
m > 2a^*.
\]
If $n$ is a positive integer and 
\beq       \label{PANT-V:defineBn}
B = \left\{  \sum_{i=0}^{n-1} a_i m^i :a_i \in A \text{ for all }  i=0,1,\ldots, n-1 \right\}
\eeq
then
\[
|B| = |A|^n
\]
\[
|B+B| = |A+A|^n
\]
\[
|B-B| = |A-A|^n.  
\]
\el

\begin{proof}
The first two identities follow immediately from the uniqueness of the $m$-adic 
representation of an integer.

If $y \in B-B$, then there exist $x = \sum_{i=0}^{n-1} a_i m^i \in B$ 
and $\tilde{x} = \sum_{i=0}^{n-1} \tilde{a}_i m^i \in B$ such that 
\[
y = x-\tilde{x}  = \sum_{i=0}^{n-1} (a_i - \tilde{a}_i ) m^i  = \sum_{i=0}^{n-1} d_i m^i  
\]
where $d_i \in A-A$ for all $i=0,1,\ldots, n-1$.

Let $d_i, d'_i  \in A-A$ for $i=0,1,\ldots, n-1$.  
We have $|d_i| \leq a^*$, $|d'_i| \leq a^*$, and so 
\[
|d_i - d'_i| \leq 2a^* \leq m-1.
\]
Define  $y, y' \in B-B$ by  
$y = \sum_{i=0}^{n-1} d_i m^i$ 
and  $y' = \sum_{i=0}^{n-1} d'_i m^i$.
Suppose that $y = y'$.  If $d_{r-1} \neq d'_{r-1}$ for some $r \in \{1,\ldots, n\}$
and $d_i = d'_i$ for $i= r, \ldots, n-1$, then
\[
0 = y-y' = \sum_{i=0}^{n-1} (d_i - d'_i ) m^i = \sum_{i=0}^{r-1} (d_i - d'_i ) m^i 
\]
and so
\[
(d'_{r-1} - d_{r-1})  m^{r-1}   = \sum_{i=0}^{r-2} (d_i - d'_i ) m^i .
\]
Taking the absolute value of each side of this equation, we obtain
\begin{align*}
m^{r-1}   & \leq |d'_{r-1} - d_{r-1}|m^{r-1} \\
&  = \left|  \sum_{i=0}^{r-2} (d_i - d'_i ) m^i \right| \\
&  \leq  2a^* \sum_{i=0}^{r-2} m^i \\
& < \left(  \frac{ 2a^*  }{m-1} \right)  m^{r-1} 
 \leq  m^{r-1}
\end{align*} 
which is absurd.  
Therefore, $y = y'$ if and only if $d_i = d'_i$ for all $i=0,1,\ldots, n-1$, 
and so $|B-B| = |A-A|^n$.
This completes the proof.
\end{proof}

Hegarty and Miller~\cite{hega-mill09} and Martin and O'Bryant~\cite{mart-obry07a} 
used probability arguments to prove that there are infinitely many MSTD sets of cardinality $k$ 
for all sufficiently large $k$.  
The following Theorem gives a constructive proof  
that, for infinitely many $k$, there exist infinitely many affinely inequivalent 
MSTD sets of integers of cardinality $k$.

\bt                   \label{PANT-V:theorem:infiniteHk}
If there exists an MSTD set of integers of cardinality $k$, then \\ 
$\mch(k^n+1) = \infty$ 
for all integers $n \geq k$.  
\et

\begin{proof}
For all integers $n \geq k \geq 1$, we have $2k-1 \geq k$ and 
\[
n(2k-1)^{n-1} \geq k\cdot k^{n-1} = k^n.
\]

Let $A$ be a nonempty set of integers of cardinality $k$.   
After an affine transformation, we can assume that $\min(A) = 0$, 
$\gcd(A) = 1$,  and $\max(A) = a^*$.  
Moreover, 
\[
A - A \supseteq \left\{ 0\} \cup \{ \pm a: a \in A \setminus \{0\} \right\}
\]
and so 
\[
|A-A| \geq 2k-1.
\]
Choose $m > 2a^*$ and $n \geq k$, 
and define the set $B$ by~\eqref{PANT-V:defineBn}.  

If $A$ is an MSTD set, then $|A+A| \geq |A-A| + 1$.  
Applying Lemma~\ref{PANT-V:lemma:Delta-m}, we obtain 
\[
|B| = k^n 
\]
and
\begin{align*}
|B+B| & = |A+A|^n \\
& \geq (|A-A|+1)^n \\
& >  |A-A|^n +  n|A-A|^{n-1} \\
& \geq  |A-A|^n +  n(2k-1)^{n-1} \\
& \geq  |A-A|^n +  k^n \\
& = |B-B| + |B|.
\end{align*}
Applying Lemma~\ref{PANT-V:lemma:Delta-k} with $\ell = k^n$, 
we see that $B$ is an MSTD set.  
Because we have infinitely many choices of $m$ and $n$, 
it follows that  $\mch(k^n+1) = \infty$.  
This completes the proof.  
\end{proof}

\bprob
Compute the smallest $k$ such that $\mch(k) = \infty$.  We know only that $k \geq 9$.  
\eprob

\bprob 
Do there exist infinitely many affinely inequivalent 
MSTD sets of integers of cardinality $k$ for all sufficiently large $k$?  
\eprob

\section{An incomplete history} 
John Marica~\cite{mari69} wrote the first paper on sets with more sums than differences.  
His paper starts with a quotation from unpublished mimeographed notes 
of  Croft~\cite{crof67}:
\begin{quotation}
\small{Problem 7 of Section VI of H. T. Croft's ``Research Problems'' (August, 1967 edition) is by J. H. Conway:

$A$ is a finite set of integers $\{ a_i \}$.  $A+A$ denotes $\{a_i + a_j\}$,   
$A - A$ denotes $\{a_i - a_j\}$.  
Prove that $A-A$ always has more numbers than $A+A$ unless $A$ is symmetrical about 0.}
\end{quotation}

I have been unable to obtain a copy of these notes.  Conway (personal communication) 
says that he did not make this conjecture, and, in fact, produced a counterexample.  
The smallest MSTD set is $ \{0, 2, 3, 4, 7, 11, 12, 14\}$, but I do not know where 
this set first appeared.  
The first published example of an MSTD set is Marica's set 
$\{1,2,3,5,8,9,13,15,16\}$.  There is a related note of Spohn~\cite{spoh71}. 
Freiman and Pigarev~\cite{frei-piga73} is another significant early work.  

Nathanson~\cite{nath07xyz} introduced the term \emph{MSTD sets}. 
There is important early work of Roesler~\cite{roes00} and 
Ruzsa~\cite{ruzs79,ruzs84,ruzs92c}, and the related paper of 
Hennecart, Robert, and Yudin~\cite{henn-robe-yudi99}.
Steve Miller and his students and colleagues have contributed greatly to this subject 
(cf.~\cite{do15a,do15b,iyer-laza-mill-zhan12, iyer-laza-mill-zhan14,mill-oros-sche10,
mill-robi-pega12, mill-sche10,zhao10a,zhao10b,zhao11}).

There has also been great interest in the Lebesgue measure of sum and difference sets  
(e.g. Steinhaus~\cite{stei20}, Piccard~\cite{picc39,picc40,picc42}, and Oxtoby~\cite{oxto71}).

\def\cprime{$'$} \def\cprime{$'$} \def\cprime{$'$} \def\cprime{$'$}
\providecommand{\bysame}{\leavevmode\hbox to3em{\hrulefill}\thinspace}
\providecommand{\MR}{\relax\ifhmode\unskip\space\fi MR }
\providecommand{\MRhref}[2]{%
  \href{http://www.ams.org/mathscinet-getitem?mr=#1}{#2}
}
\providecommand{\href}[2]{#2}

\end{document}